\documentclass[12pt,reqno]{article}

\usepackage[usenames]{color}
\usepackage{amssymb}
\usepackage{amsmath}
\usepackage{amsthm}
\usepackage{amsfonts}
\usepackage{graphicx}
\usepackage{tikz}
\usepackage{mathtools}
\usetikzlibrary{decorations.markings}

\usepackage{hyperref}

\usepackage{color}
\usepackage{fullpage}
\usepackage{float}

\usepackage{graphicx}
\usepackage{latexsym}
\usepackage{epsf}
\usepackage{breakurl}

\usepackage{tikz}
\usetikzlibrary{calc}

\begin{document}
	
	\theoremstyle{plain}
	\newtheorem{theorem}{Theorem}
	\newtheorem{corollary}[theorem]{Corollary}
	\newtheorem{lemma}[theorem]{Lemma}
	\newtheorem{proposition}[theorem]{Proposition}
	
	\theoremstyle{definition}
	\newtheorem{definition}[theorem]{Definition}
	\newtheorem{example}[theorem]{Example}
	\newtheorem{conjecture}[theorem]{Conjecture}
	\newtheorem{notation}[theorem]{Notation}
	\newtheorem{remark}[theorem]{Remark}
	\newtheorem{claim}[theorem]{Claim}
	
	\theoremstyle{remark}

	\begin{center}
		\vskip 1cm{\Large\bf 
			Cyclic and Linear Graph Partitions \\ and Normal Ordering
		}
		\vskip 1cm
		\large
		Ken Joffaniel Gonzales\\
		Department of Physical Sciences and Mathematics\\
		University of the Philippines Manila\\
		Manila, Philippines\\
		\href{kmgonzales1@up.edu.ph}{\tt kmgonzales1@up.edu.ph} \\
	\end{center}
	
	\vskip .2in
	
	\begin{abstract}
		The Stirling number of a simple graph is the number of partitions of its vertex set into a specific number of non-empty independent sets. In 2015, Engbers et al.\ showed that the coefficients in the normal ordering of a word $w$ in the alphabet $\{x,D\}$ subject to the relation $Dx=xD+1$ are equal to the Stirling number of certain graphs constructed from $w$. In this paper, we introduce graphical versions of the Stirling numbers of the first kind and the Lah numbers and show how they occur as coefficients in other normal ordering settings. Identities involving their $q$-analogues are also obtained.
	\end{abstract}
	
	\section{Introduction}\label{sec:intro}
	
	The \emph{Weyl algebra} is the noncommutative algebra generated by operators $x$ and $D$ satisfying $x(f(x))=xf(x)$ and $D(f(x))=\frac{d}{dx}f(x)$. This algebra is also defined via the relation $Dx=xD+1$ without reference to any specific representation of $x$ and $D$ as operators. The \emph{normally ordered form} of a word $w$ in the Weyl algebra with $m$ $x$'s and $n$ $D$'s is the unique expansion given by 
	\[
	w = x^{m-n}\sum_{k=0}^n {w \brace k} x^k D^k\,.
	\]
	
	The coefficients $\displaystyle {w \brace k}$ have a well-known interpretation in terms of rook placements~\cite{Nav,Var}. More recently, Engbers et al.~\cite{Engetal} derived another interpretation by showing that for some graph $G$ constructed from $w$, the number of partitions of $V(G)$ into $k$ non-empty independent sets are exactly the coefficients $\displaystyle {w \brace k}$. If $w=(xD)^n$, the associated graph is the empty graph with $n$ vertices and the partitions coincide with the usual interpretation for the Stirling numbers of the second kind $S(n,k)$. 
	
	This paper seeks to extend the results of Engbers et al.~\cite{Engetal} by introducing graphical versions of the Stirling numbers of the first kind $c(n,k)$ and the Lah numbers $L(n,k)$ that are inspired by their classical interpretations. We then use these numbers to combinatorially interpret certain normal order coefficients.
	
	The paper is organized as follows. In Section~\ref{sec:revisit}, we recall the constructions and main results of~\cite{Engetal}. We then introduce rook placements associated with graphs in Section~\ref{sec:rook}. The graphical versions of $c(n,k)$ and $L(n,k)$ involve cyclic and linear partitions, and these are discussed in Sections~\ref{sec:cycle} and~\ref{sec:line}, respectively.
	
	\section{Graphical Stirling numbers of the second kind} \label{sec:revisit}
	
	We begin by recalling some terminology and results from~\cite{Engetal} relating the coefficients $\displaystyle {w \brace k}$ and vertex partitions of some graphs constructed from $w$.
	
	\begin{definition}
		A word $w$ is a \emph{Dyck word} if the number of $x$'s is equal to the number of $D$'s and for every $k$, every subword consisting of the first to the $k$-th letter from the left of $w$ has a least as many $x$'s as $D$'s. A Dyck word is \emph{irreducible} if either $w=xD$ or $w=xw'D$ for some non-empty Dyck word. Otherwise, the Dyck word is \emph{reducible}.
	\end{definition}
	
	Given a Dyck word $w$, we associate an unlabeled simple graph $G_w$ constructed as follows.
	\begin{enumerate}
		\item If $w=x D$, then $G_w$ is a graph with one vertex. 
		\item If $w$ is irreducible and $w=xw'D$ for some non-empty Dyck word $w'$, then $G_w$ is the graph obtained by adding a single dominating vertex to $G_{w'}$.
		\item If $w$ is reducible, say $w=w_1w_2\cdots w_m$ where each $w_{i}$ is irreducible, then $G_w$ is the disjoint union of the $G_{w_i}$'s.
	\end{enumerate}
	
\begin{figure}[ht!] 
		\centering
		\begin{tikzpicture}[scale=0.8, transform shape]
			\begin{scope}[yscale=1,xscale=-1,xshift=-12cm]
				\node (n1) at (0,0) [circle,draw=black] {};
				\node (n2) at (2,0) [circle,draw=black] {};
				\node (n3) at (4,0) [circle,draw=black] {};
				\node (n4) at (6,0) [circle,draw=black] {};
				\node (n5) at (8,0) [circle,draw=black] {};
				\node (n6) at (10,0) [circle,draw=black] {};
				\node (n7) at (12,0) [circle,draw=black] {};

				\draw (n1) .. controls (1,1.5) and (3,1.5) .. (n3);
				\draw (n2) -- (n3);
				\draw (n1) .. controls (2,3.5) and (6,3.5) .. (n5);
				\draw (n2) .. controls (3.4,2.3) and (6.6,2.3) .. (n5);
				\draw (n3) .. controls (5,1.5) and (7,1.5) .. (n5);
				\draw (n4) -- (n5);
				\draw (n6) -- (n7);
			\end{scope}
			
			\node at (0,-0.75) {\small 1};
			\node at (2,-0.75) {\small 2};
			\node at (4,-0.75) {\small 3};
			\node at (6,-0.75) {\small 4};
			\node at (8,-0.75) {\small 5};
			\node at (10,-0.75) {\small 6};
			\node at (12,-0.75) {\small 7};
			
		\end{tikzpicture}
		\caption{The graph $G_w$ associated with $w=xxxDxDDxDDxxDD$} \label{fig:gw}
	\end{figure}
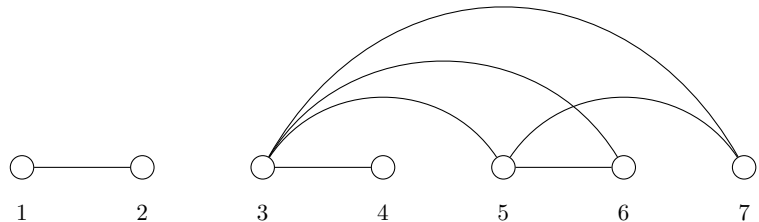

	Figure~\ref{fig:gw} shows the graph $G_w$ where $w=xxxDxDDDxxDD$. While we take $G_w$ to be an unlabeled graph, the figure shows a labeling of the vertices that will be relevant in Section~\ref{sec:rook}.
	
	\begin{definition}
		The $k$-th \emph{Stirling number} of a simple graph $G$, denoted by $S(G,k)$ is the number of ways to partition $V(G)$ into $k$ non-empty independent sets.
	\end{definition}
	
	Denote the length of a word $w$ by $|w|$. Using the representation of the Weyl algebra in terms of differential operators, Engbers et al.~\cite{Engetal} showed that if $w$ is a Dyck word with $|w|/2=n$, then
	\[
	w = \sum_{k=0}^n S(G_w,k)x^k D^k\,.
	\]
	That is, $\displaystyle {w\brace k}= S(G_w,k)$. The authors also obtained an interpretation for the normal order coefficients of arbitrary words by adjoining appropriate powers of $x$ and $D$ to form a Dyck word~\cite[Corollary 2.5]{Engetal}. As we shall be eventually defining the first kind version of $S(G,k)$, we shall call the number $S(G,k)$ as $k$-th \emph{Stirling number of the second kind of $G$} from here onwards.
	
\begin{figure}[ht!] 
		\centering
		\begin{tikzpicture}[scale=0.8, transform shape]
			\filldraw[color=gray!50, fill=gray!50] (0,3) rectangle (1,7);
			\filldraw[color=gray!50, fill=gray!50] (1,4) rectangle (3,7);
			\filldraw[color=gray!50, fill=gray!50] (3,5) rectangle (5,7);
			
			\foreach \i in {0,...,7}{
				\draw (0,\i) -- (7,\i);
				\draw (\i,0) -- (\i,7);
			}
			
			\draw[ultra thick] (0.025,0) -- (0.025,3) -- (1,3) -- (1,4) -- (3,4) -- (3,5) -- (5,5) -- (5,6) -- (5,6.975) -- (7,6.975);
			
			\draw[dotted] (0,0) -- (7,7);
			
			\foreach \i in {1,...,7}{
				\node at (-0.3,7.5-\i) {\scriptsize $\i$};
				\node at (7.5-\i,7.3) {\scriptsize $\i$};
			}
			
			\draw[dashed] (0,0) -- (7,7);
			
			\node at (0.5,1.5) {\scriptsize $(7,6)$};
			\node at (0.5,2.5) {\scriptsize $(7,5)$};
			\node at (1.5,2.5) {\scriptsize $(6,5)$};
			\node at (1.5,3.5) {\scriptsize $(6,4)$};
			\node at (2.5,3.5) {\scriptsize $(5,4)$};
			\node at (3.5,4.5) {\scriptsize $(4,3)$};
			\node at (5.5,6.5) {\scriptsize $(2,1)$};
		\end{tikzpicture}
		\caption{The Dyck path $\mathcal P(w)$ (thick lines) and board $B(w)$ (gray cells) associated with {$w=xxxDxDDxDDxxDD$.}} \label{fig:bw}
		
	\end{figure}
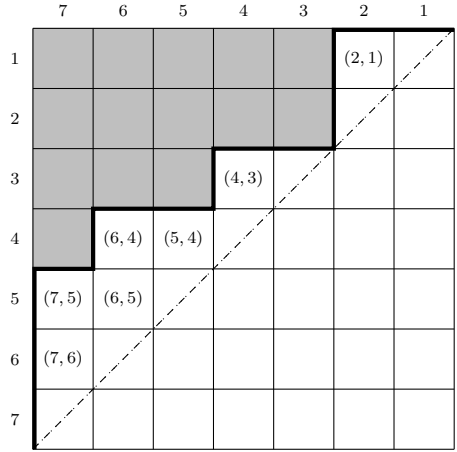

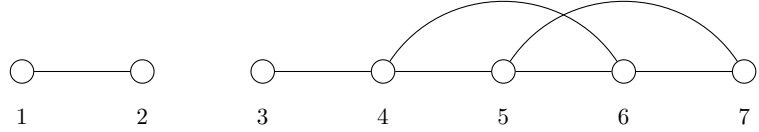
\begin{figure}[ht!] 
		\centering
		\begin{tikzpicture}[scale=0.8, transform shape]
			\begin{scope}[yscale=1,xscale=-1,xshift=-12cm]
				\node (n1) at (0,0) [circle,draw=black] {};
				\node (n2) at (2,0) [circle,draw=black] {};
				\node (n3) at (4,0) [circle,draw=black] {};
				\node (n4) at (6,0) [circle,draw=black] {};
				\node (n5) at (8,0) [circle,draw=black] {};
				\node (n6) at (10,0) [circle,draw=black] {};
				\node (n7) at (12,0) [circle,draw=black] {};
				
				\draw (n1) -- (n2);
				\draw (n2) -- (n3);
				\draw (n3) -- (n4);
				\draw (n4) -- (n5);
				\draw (n6) -- (n7);
				\draw (n1) .. controls (1,1.5) and (3,1.5) .. (n3);			
				\draw (n2) .. controls (3,1.5) and (5,1.5) .. (n4);
			\end{scope}
			
			\node at (0,-0.75) {\small 1};
			\node at (2,-0.75) {\small 2};
			\node at (4,-0.75) {\small 3};
			\node at (6,-0.75) {\small 4};
			\node at (8,-0.75) {\small 5};
			\node at (10,-0.75) {\small 6};
			\node at (12,-0.75) {\small 7};
		\end{tikzpicture}
		\caption{The graph $H_w$ associated with $w=xxxDxDDxDDxxDD$} \label{fig:hw}
		
	\end{figure}

	In addition, there is another graph $H_w$ for which $\displaystyle {w\brace k}=S(H_w,k)$. To construct $H_w$, we first draw the \emph{Dyck path} $\mathcal P(w)$ starting from the lower left corner of the $n \times n$ lattice, where $n=|w|/2$. The steps on $\mathcal P(w)$ are then obtained by replacing $x$ with a unit up step and $D$ with a unit right step. Observe that since $w$ is a Dyck word, its Dyck path $\mathcal P(w)$ ends in the upper right corner and may touch but not cross the northeast diagonal of the lattice. Next, we associate with $w$ the board $B(w)$ consisting of cells in the lattice lying above $\mathcal P(w)$. Label the columns of the lattice with $1,2,\ldots,n$ from right to left and its rows $1,2,\ldots,n$ from top to bottom and denote the cell in column $i$ and row $j$ by $(i,j)$. We construct an edge from $i$ to $j$ in $H_w$ for each $i>j$ if the cell $(i,j)$ lies outside $B(w)$ and above the northeast diagonal. Figure~\ref{fig:bw} shows the Dyck path and board with $w=xxxDxDDxDDxxDD$ while Figure~\ref{fig:hw} shows the corresponding graph $H_w$. (Note that we will be using this $w$ in all the examples in the succeeding sections.)
	
	Engbers et al.~\cite{Engetal} showed two proofs that $\displaystyle {w\brace k}=S(H_w,k)$. The first proof uses results involving rook placements by Navon~\cite{Nav}, Varvak~\cite{Var} and Goldman et al.~\cite{Goletal}. The second proof is direct and proceeds by showing that $G_w$ and $H_w$ have the same chromatic polynomials.
	
	\section{Boards associated with graphs} \label{sec:rook}
	
	In this section, we introduce boards associated with graphs which will be used in the succeeding sections. Additionally, they provide a visualization of why the number of partitions of $H_w$ and $G_w$ both produce the normal order coefficients.
	
	\begin{definition} A \emph{placement of non-attacking rooks} on a board $B$ is a marking of the cells of $B$ with rooks ``$\times$'' such that no two rooks lie in the same row or column. The $k$-th \emph{rook number} of $B$, denoted by $r_k(B)$, is the number of placements of $k$ non-attacking rooks on $B$. Furthermore, we say that two boards $B_1$ and $B_2$ are \emph{rook equivalent} (or simply, \emph{equivalent}) if $r_k(B_1)=r_k(B_2)$ for every $k$.
	\end{definition}
	
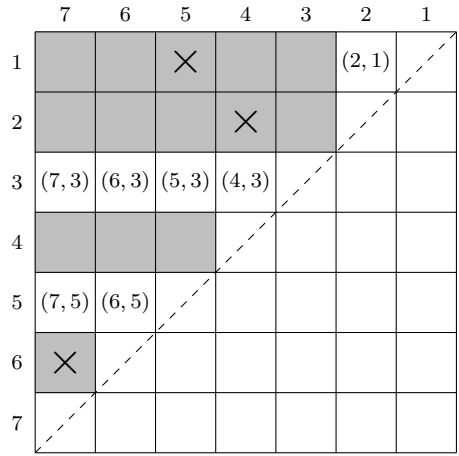
\begin{figure}[ht!] 
		\centering
		\begin{tikzpicture}[scale=0.8]
			\filldraw[color=gray!50, fill=gray!50] (0,1) rectangle (1,2);
			\filldraw[color=gray!50, fill=gray!50] (0,3) rectangle (3,4);
			\filldraw[color=gray!50, fill=gray!50] (0,5) rectangle (5,7);
			
			\foreach \i in {0,...,7}{
				\draw (0,\i) -- (7,\i);
				\draw (\i,0) -- (\i,7);
			}

			\foreach \i in {1,...,7}{
				\node at (-0.3,7.5-\i) {\scriptsize $\i$};
				\node at (7.5-\i,7.3) {\scriptsize $\i$};
			}

			\draw[dashed] (0,0) -- (7,7);
			
			\node at (0.5,1.5) {\Large $\times$};
			\node at (2.5,6.5) {\Large $\times$};
			\node at (3.5,5.5) {\Large $\times$};
			
			\node at (0.5,2.5) {\scriptsize $(7,5)$};
			\node at (1.5,2.5) {\scriptsize $(6,5)$};
			\node at (0.5,4.5) {\scriptsize $(7,3)$};
			\node at (1.5,4.5) {\scriptsize $(6,3)$};
			\node at (2.5,4.5) {\scriptsize $(5,3)$};
			\node at (3.5,4.5) {\scriptsize $(4,3)$};
			\node at (5.5,6.5) {\scriptsize $(2,1)$};				
		\end{tikzpicture} 
		\caption{A placement of $3$ non-attacking rooks on a board.} \label{fig:rookboard}
		
	\end{figure}

	Figure~\ref{fig:rookboard} shows a placement of $3$ non-attacking rooks.
	
	\begin{definition} Given a simple labeled graph $G$ with vertices labeled $1,2,\ldots,|G|$, the \emph{board associated with $G$}, denoted by $\mathcal B(G)$, is the board consisting of cells $(i,j)$ if vertices $i$ and $j$ are not adjacent in $G$ and $i>j$.
	\end{definition}	
	
	By construction, the vertices of $H_w$ are labeled, with the edges added later as determined by $B(w)$. We take this labeling as the prescribed labeling for $H_w$. For $G_w$, we define its prescribed labeling as follows. Recall that a pair of $x$ and $D$ in $w$ results into either an additional isolated vertex or an additional dominating vertex so that we can associate each $D$ in $w$ with its corresponding vertex in $G_w$. Label the $D$'s in $w$ from right to left by $1,2,\ldots,|w|/2$. The label of the vertices of $G_w$ are then the label of the corresponding $D$. Figure~\ref{fig:gw} shows a graph $G_w$ in this prescribed labeling.
	
	For a simple graph $G$, the value of $S(G,k)$ is independent of its labeling. However, two different labelings of $G$ do not always the result into equivalent boards. We therefore assume from here onwards that $H_w$ and $G_w$ are labeled in the prescribed manner. By construction of $H_w$, we have that $\mathcal B(H_w)=B(w)$, a property not satisfied by $G_w$ in general. Figure~\ref{fig:rookboard}, minus the rooks, shows $\mathcal B(G_w)$ for $w=xxxDxDDxDDxxDD$. Observe that flushing the rows of $\mathcal B(G_w)$ in the figure towards the top produces $\mathcal B(H_w)$ and thus, the two boards are equivalent. We prove that this is true for any Dyck word $w$ in the theorem that follows.
	
	\begin{theorem} \label{thm:alt} Let $w$ be a Dyck word. Then,
		\begin{enumerate}
			\item For every $i$, column $i$ of $\mathcal B(H_w)$ and $\mathcal B(G_w)$ have the same number of cells.
			\item For every $j$, row $j$ of $\mathcal B(H_w)$ and the $j$-th non-empty row from the top of $\mathcal B(G_w)$ have the same number of cells.
			\item The boards $\mathcal B(H_w)$ and $\mathcal B(G_w)$ are equivalent.
		\end{enumerate}
		
		\begin{proof} We begin by looking at $\mathcal B(H_w)=B(w)$. The word $w$ determines the shape of $B(w)$, and hence, the number of cells in each column. Specifically, the number of cells in column $i$ is equal to the number $x$'s to the right of the $i$-th $D$ from the right.
			
			Next, we look at $\mathcal B(G_w)$. If $w$ is a Dyck word, then the $x$'s and $D$'s come in pairs, with each pair resulting into either a non-dominating vertex if it occurs as $xD$ or a dominating vertex if it occurs as $xw'D$ for some non-empty word $w'$. As before, label the $D$'s $1,2,\ldots,|w|/2$ from right to left and label the $x$'s with the label of the $D$ it is paired with. Then, for $i<j$, vertex $i$ dominates (and is therefore adjacent to) vertex $j$ if and only if the $x$ and $D$ labeled $i$ sandwiches the $x$ and $D$ labeled $j$. This implies that if the $D$ labeled $i$ occurs to the left of the $x$ labeled $j$, then the $x$ and $D$ labeled $i$ do not sandwich the $x$ and $D$ labeled $j$, so that vertices $i$ and $j$ are non-adjacent. Hence, the number of $x$'s to the right of the $D$ labeled $i$ counts the number of vertices not adjacent to vertex $i$ and therefore, the number of cells in column $i$ of $\mathcal B(G_w)$. This is exactly the number of cells of column $i$ of $\mathcal B(H_w)$, thereby proving the first statement.
			
			Given two boards $B_1$ and $B_2$, it is easy to see that if the cells of both boards are flushed to the left, the second statement follows from the first statement. Since $\mathcal B(H_w)=B(w)$, the cells of $\mathcal B(H_w)$ are flushed to the left. As to $\mathcal B(G_w)$, suppose that for some $j<i_1<i_2$, $\mathcal B(G_w)$ contains the cell $(i_1,j)$ but not the cell $(i_2,j)$. This means that vertex $j$ dominates vertex $i_2$ but not $i_1$. In terms of $w$, this implies that the $D$ labeled $j$ occurs to the left of the $D$ labeled $i_1$ and to the right of the $D$ labeled $i_2$, which cannot happen since $j<i_1 < i_2$.
			
			Finally, the third statement follows from the first two.
		\end{proof}
	\end{theorem}
	
	\begin{corollary} \label{cor:alt} Let $w$ be a Dyck word. Then for any $k$, $S(G_w,k)=S(H_w,k)$.
		\begin{proof}
			If $w$ is a Dyck word with $|w|/2=n$, then a bijection between the partitions of the vertex set of $G_w$ into $k$ non-empty independent sets and the placements of $n-k$ non-attacking rooks on $\mathcal B(G_w)$ is as follows. Two elements $i$ and $j$ belong to the same subset in the partition if and only if there is a rook in cell $(i,j)$. This shows that $S(G_w,k)=r_{n-k}(\mathcal B(G_w))$. For example, the rook placement associated with the partition $\{1,5\}\{2,4\}\{3\}\{6,7\}$ of the graph $G_w$ in Figure~\ref{fig:gw} is shown in Figure~\ref{fig:rookboard}. (This bijection was originally defined in~\cite[Section 1.3]{Butetal} to be between the placement of $n-k$ rooks in the staircase board $B((xD)^n)$ and the partitions of the set $\{1,2,\ldots,n\}$ into $k$ subsets.)
			
			The same rule gives a bijection between the partitions of the vertex set of $H_w$ into $k$ non-empty independent sets and the placements of $n-k$ non-attacking rooks on $\mathcal B(H_w)$. Hence, $S(H_w,k)=r_{n-k}(\mathcal B(H_w))$. The result then follows since $\mathcal B(G_w)$ and  $\mathcal B(H_w)$ are equivalent boards by the third statement in Theorem~\ref{thm:alt}.
		\end{proof}
	\end{corollary}

	\section{Cyclic partitions} \label{sec:cycle}
	
	In this section, we consider words in the alphabet $\{V,U\}$ subject to the relation $UV=VU+V$. We write the normally ordered form of a word $w$ with $m$ $V$'s and $n$ $U$'s by
	\[
	w = V^{m}\sum_{k=0}^n {w\brack k} U^k\,.
	\]
	Both the relations $UV=VU+V$ and $Dx=xD+1$ are special cases of the relation $UV=VU+hV^s$ studied by Mansour et al.~\cite{Manetal}, where the authors also showed that the normally ordered form of $(VU)^n$ is given by
	\[
	(VU)^n = V^n \sum_{k=0}^n c(n,k) U^k\,.
	\]
	The coefficients $\begin{bmatrix} w \\ k \end{bmatrix}$ may therefore be viewed as a kind of generalization of the Stirling numbers of the first kind.
	
	A $q$-analogue of the relation $UV=VU+V$ is given by $UV=qVU+V$. We write the normally ordered form of a word $w$ with $m$ $V$'s and $n$ $U$'s under this relation by
	\[
	w = V^{m}\sum_{k=0}^n {w;q\brack k} U^k\,.
	\]
	
	Given a board $B$, we denote by $f_{k}(B)$ the number of placements of $k$ rooks on $B$ such that no two rooks lie in the same column. We call rook placements with this restriction \emph{file placements}. For a file placement $\phi$ of rooks, the \emph{weight} of $\phi$, denoted by $\mathrm {wt}(\phi)$, is the number of cells not lying above a rook. The $k$-th $q$-\emph{file number} of $B$ is then given by
	\[
	f_{n;q}(B) = \sum_{\phi} q^{\mathrm{wt}(\phi)}\,,
	\]
	where the sum runs over all file placements $\phi$ of $n$ rooks on $B$. 
	
	If a word $w$ has $n$ $U$'s, a result by Celeste et al.~\cite[Identity 3]{Celetal}, in terms of the notation in this paper, implies that
	\begin{align}
		{w\brack k} &= f_{n-k}(B(w)) \label{eqn:cel} \\
		{w;q\brack k} &= f_{n-k;q}(B(w))\,. \label{eqn:celq}
	\end{align}

	We now define our graphical version of the Stirling number of the first kind.
	
	\begin{definition} Given a permutation $\sigma$, let $\sigma^{(i)}$ be the partial permutation obtained by keeping only elements $j\leq i$ in the cycle structure of $\sigma$. A $k$-\emph{cyclic partition} of a labeled graph $G$ is a permutation $\sigma$ of $V(G)$ into $k$ cycles such that for every vertex $i$, $\sigma^{(i)}(i)=j$ only if there is no directed edge from $i$ to $j$. The \emph{$k$-th Stirling number of the first kind of $G$}, denoted by $c(G,k)$, is number of its $k$-cyclic partitions.
	\end{definition}
	
	If $G$ is the empty graph with $n$ vertices, then $c(G,k)=c(n,k)$. 
	
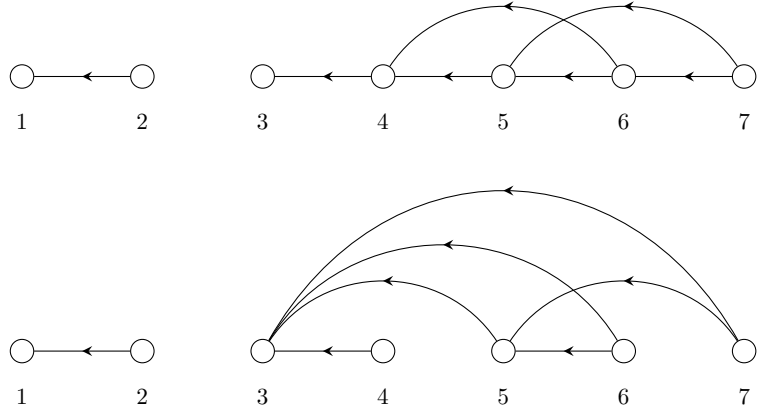
\begin{figure}[ht!] 
		\centering
		\begin{tikzpicture} [scale=0.8, transform shape, decoration={markings,mark= at position 0.5 with {\pgftransformscale{1.2}\arrow{stealth}}}]
			
			\begin{scope}[yscale=1,xscale=-1,xshift=-12cm]
				\node (n1) at (0,0) [circle,draw=black] {};
				\node (n2) at (2,0) [circle,draw=black] {};
				\node (n3) at (4,0) [circle,draw=black] {};
				\node (n4) at (6,0) [circle,draw=black] {};
				\node (n5) at (8,0) [circle,draw=black] {};
				\node (n6) at (10,0) [circle,draw=black] {};
				\node (n7) at (12,0) [circle,draw=black] {};
				
				\draw [postaction={decorate}] (n1) -- (n2);
				\draw [postaction={decorate}] (n2) -- (n3);
				\draw [postaction={decorate}] (n3) -- (n4);
				\draw [postaction={decorate}] (n4) -- (n5);
				\draw [postaction={decorate}] (n6) -- (n7);
				\draw [postaction={decorate}] (n1) .. controls (1,1.5) and (3,1.5) .. (n3);			
				\draw [postaction={decorate}] (n2) .. controls (3,1.5) and (5,1.5) .. (n4);			
			\end{scope}	
			
			\node at (0,-0.75) {\small 1};
			\node at (2,-0.75) {\small 2};
			\node at (4,-0.75) {\small 3};
			\node at (6,-0.75) {\small 4};
			\node at (8,-0.75) {\small 5};
			\node at (10,-0.75) {\small 6};
			\node at (12,-0.75) {\small 7};
		\end{tikzpicture}		
		
		\begin{tikzpicture} [scale=0.8, transform shape, decoration={markings,mark= at position 0.5 with {\pgftransformscale{1.2}\arrow{stealth}}}]
			
			\begin{scope}[yscale=1,xscale=-1,xshift=-12cm]	
				\node (n1) at (0,0) [circle,draw=black] {};
				\node (n2) at (2,0) [circle,draw=black] {};
				\node (n3) at (4,0) [circle,draw=black] {};
				\node (n4) at (6,0) [circle,draw=black] {};
				\node (n5) at (8,0) [circle,draw=black] {};
				\node (n6) at (10,0) [circle,draw=black] {};
				\node (n7) at (12,0) [circle,draw=black] {};
				
				\draw [postaction={decorate}] (n1) .. controls (1,1.5) and (3,1.5) .. (n3);
				\draw [postaction={decorate}] (n2) -- (n3);
				\draw [postaction={decorate}] (n1) .. controls (2,3.5) and (6,3.5) .. (n5);
				\draw [postaction={decorate}] (n2) .. controls (3.4,2.3) and (6.6,2.3) .. (n5);
				\draw [postaction={decorate}] (n3) .. controls (5,1.5) and (7,1.5) .. (n5);
				\draw [postaction={decorate}] (n4) -- (n5);
				\draw [postaction={decorate}] (n6) -- (n7);
			\end{scope}
			
			\node at (0,-0.75) {\small 1};
			\node at (2,-0.75) {\small 2};
			\node at (4,-0.75) {\small 3};
			\node at (6,-0.75) {\small 4};
			\node at (8,-0.75) {\small 5};
			\node at (10,-0.75) {\small 6};
			\node at (12,-0.75) {\small 7};
		\end{tikzpicture}
		
		\caption{The labeled directed graphs $H^d_w$ (top) and 	$G^d_w$ (bottom) associated with {$w=xxxDxDDxDDxxDD$.}} \label{fig:dir}
		
	\end{figure}

	In order to state our interpretation for $\begin{bmatrix} w \\ k \end{bmatrix}$, we extend the definition of a Dyck word in $\{V,U\}$, as well as the graphs $H_w$ and $G_w$ by adapting their definitions in terms of $\{x,D\}$, with $V,U$ replacing $x,D$, respectively. For a Dyck word $w$, we denote by $H^d_w$ the labeled digraph obtained from $H_w$ by replacing every edge connecting $i$ to $j$, $i<j$, with a directed edge from $j$ to $i$. The graph $G^d_w$ is defined similarly. The directed versions of the graphs in Figures~\ref{fig:hw}~and~\ref{fig:gw} are shown in Figure~\ref{fig:dir}.
	
	We also introduce the following map, which is a slight modification of a map in \cite[Section 1.3]{Butetal}. Let $G$ be either $H^d_w$ or $G^d_w$ for some Dyck word $w$ such that $|w|/2=n$. Given a cyclic partition $\sigma$ of $G$, let $\Phi(\sigma)$ be the file placement of rooks on $\mathcal B(G)$ obtained as follows. For $i=1,2,\ldots,n$:
	\begin{enumerate}
		\item If $i$ is a minimal element, place no rook in column $i$ of $\mathcal B(G)$ .
		\item If $i$ is not a minimal element, then place a rook in cell $(i,\sigma^{(i)}(i))$.
	\end{enumerate}
	
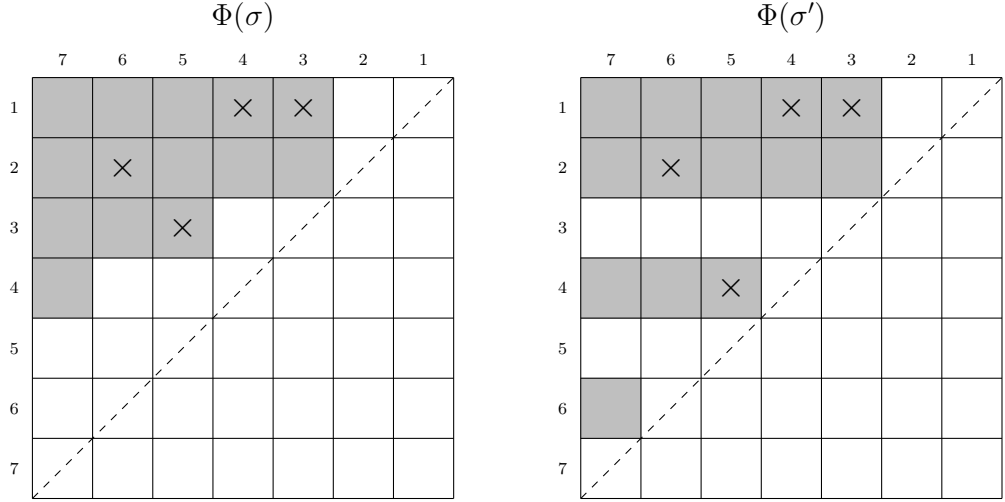
\begin{figure}[ht!] 
		\centering
		\begin{tikzpicture}[scale=0.8, transform shape]
			
			\filldraw[color=gray!50, fill=gray!50] (0,3) rectangle (1,7);
			\filldraw[color=gray!50, fill=gray!50] (1,4) rectangle (3,7);
			\filldraw[color=gray!50, fill=gray!50] (3,5) rectangle (5,7);
			
			\foreach \i in {0,...,7}{
				\draw (0,\i) -- (7,\i);
				\draw (\i,0) -- (\i,7);
			}
			
			\foreach \i in {1,...,7}{
				\node at (-0.3,7.5-\i) {\scriptsize $\i$};
				\node at (7.5-\i,7.3) {\scriptsize $\i$};
			}
			
			\draw[dashed] (0,0) -- (7,7);
			
			\node at (1.5,5.5) {\Large $\times$};
			\node at (2.5,4.5) {\Large $\times$};
			\node at (3.5,6.5) {\Large $\times$};
			\node at (4.5,6.5) {\Large $\times$};
			
			\node at (3.5,8) {\large $\Phi(\sigma)$};
		\end{tikzpicture}
		\hspace{1cm}
		\begin{tikzpicture}[scale=0.8, transform shape]
			\filldraw[color=gray!50, fill=gray!50] (0,1) rectangle (1,2);
			\filldraw[color=gray!50, fill=gray!50] (0,3) rectangle (3,4);
			\filldraw[color=gray!50, fill=gray!50] (0,5) rectangle (5,7);
			
			\foreach \i in {0,...,7}{
				\draw (0,\i) -- (7,\i);
				\draw (\i,0) -- (\i,7);
			}
			
			\foreach \i in {1,...,7}{
				\node at (-0.3,7.5-\i) {\scriptsize $\i$};
				\node at (7.5-\i,7.3) {\scriptsize $\i$};
			}
			
			\draw[dashed] (0,0) -- (7,7);
			
			\node at (1.5,5.5) {\Large $\times$};
			\node at (2.5,3.5) {\Large $\times$};
			\node at (3.5,6.5) {\Large $\times$};
			\node at (4.5,6.5) {\Large $\times$};
			
			\node at (3.5,8) {\large $\Phi(\sigma')$};
		\end{tikzpicture} 
		\caption{The file placement of rooks $\Phi(\sigma)$ associated with the cyclic partition $\sigma=(1534)(26)(7)$ of $H^d_w$ and the file placement of rooks $\Phi(\sigma')$ associated with the cyclic partition $\sigma'=(1354)(26)(7)$ of $G^d_w$, with {$w=xxxDxDDxDDxxDD$.}} \label{fig:phi}
		
	\end{figure}

	As illustration, consider the graphs $H^d_w$ and $G^d_w$ shown in Figure~\ref{fig:dir}. A $3$-cyclic partition of $H^d_w$ is given by $\sigma=(1534)(26)(7)$. Columns $1$, $2$ and $7$ of $\mathcal B(H_w)$ do not contain a rook since these are the minimal elements of $\sigma$. In addition, $\sigma^{(3)}(3)=1$, $\sigma^{(4)}(4)=1$, $\sigma^{(5)}(5)=3$ and $\sigma^{(6)}(6)=2$, so that the associated rook placement has a rook in cells $(3,1)$, $(4,1)$, $(5,3)$ and $(6,2)$. On the other hand, a $3$-cyclic partition of $G^d_w$ is given by $\sigma'=(1354)(26)(7)$. The rook placements associated with $\sigma$ and $\sigma'$ are shown in Figure~\ref{fig:phi}.
	
	We can now state our main result for this section. 
	
	\begin{theorem} \label{thm:firstmain} Let $w$ be a Dyck word in $\{V,U\}$ subject to $UV=VU+V$. Then $\begin{bmatrix} w \\ k \end{bmatrix} = c(H^d_w,k) = c(G^d_w,k)$.
		
		\begin{proof} First, we show that $\Phi$ is bijection between $k$-cyclic partitions of $H^d_w$ with $|w|/2=n$ and the file placement of $n-k$ rooks on $\mathcal B(H_w)=B(w)$. The result $\begin{bmatrix} w \\ k \end{bmatrix} = c(H^d_w,k)$ then follows from Equation~\eqref{eqn:cel}. 
			
			Let $\sigma$ be a $k$-cyclic partition of $H^d_w$. Suppose $i$ is not a minimal element of $\sigma$ and that column $i$ of $B(w)$ has $j$ cells. Then, there is an edge from $i$ to each of $i-1,i-2,\ldots,j+1$ and $\sigma^{(i)}(i)$ cannot be equal to any of these by the definition of a cyclic partition. Instead, $\sigma^{(i)}(i)$ is equal to one of $1,2\ldots,j$, each of which corresponds to a placement of a rook in one of these $j$ cells, specifically at cell $(i,\sigma^{(i)}(i))$, $\sigma^{(i)}(i)=1,2,\ldots,j$.  Moreover, since $\sigma$ has $k$ minimal elements, $\Phi(\sigma)$ contains $k$ columns without rooks, so that $\Phi(\sigma)$ contains $n-k$ rooks.
			
			Next, we establish $c(H_w,k) = c(G_w,k)$ by giving a bijection  between the $k$-cyclic partitions of $H^d_w$ and the $k$-cyclic partitions of $G^d_w$. List the rows of $\mathcal B(G_w)$ by $t_1 < t_2 < \ldots < t_{m_1}, t_{m_1+1} < t_{m_1+2} < \ldots < t_n$, where rows $t_1,t_2,\ldots,t_{m_1}$ have nonzero length and rows $t_{m_1+1}, t_{m_1+2}, \ldots,t_n$ have zero length. Let $\sigma$ be a $k$-cyclic partition of $H^d_w$ and $\sigma'$ be the cyclic partition having the same minimal elements as $\sigma$ and satisfying  $\sigma'^{(i)}(i)=t_{\sigma^{(i)}(i)}$ for every non-minimal element $i$. This construction combined with the second statement of Theorem~\ref{thm:alt} shows that the rook placements associated with $\sigma$ and $\sigma'$ are identical up to a reordering of the rows. More precisely, for every $j$, the configuration of row $j$ of $\mathcal B(H_w)$ is the same as the configuration of row $t_j$ of $\mathcal B(G_w)$. Thus, $\sigma'$ is a $k$-cyclic partition of $G^d_w$ and priming gives the desired bijection. For example, the listing of the rows of the board $\mathcal B(G_w)$ shown in Figure~\ref{fig:rookboard} in the order just described is given by $1<2<4<6,3<5<7$. If $\sigma=(1534)(26)(7)$, then $\sigma'=(1354)(26)(7)$.
		\end{proof}
	\end{theorem}
	
	The next theorem gives a recurrence formula for $c(H^d_w,k)$.
	
	\begin{corollary}\label{cor:rec1} \label{thm:cnkrec} Let $w$ be a Dyck word in $\{V,U\}$ such that $n=|w|/2$, subject to $UV=VU+V$. Suppose that the column lengths of $B(w)$ are $b_1,b_2,\ldots,b_n$. Let $w_i$ be a subword of $w$ such that the board $B(w_i)$ consists of columns $1,2,\ldots,i$ of $B(w)$. Then for $k\leq m\leq n$,
		\[
		c(H^d_{w_m},k) = c(H^d_{w_{m-1}},k-1) + b_{m} c(H^d_{w_{m-1}},k)\,.
		\]	
		
		\begin{proof} By Theorem~\ref{thm:firstmain} and Equation~\eqref{eqn:cel},
			\[
			f_{m-k}(B(w))=c(H^d_{w_m},k)\,.
			\]
			The number of file placements of $m-k$ rooks on $B(w_m)$ without a rook in the $m$-th column is equal to the number of file placements of $m-k$ rooks on $B(w_{m-1})$, which is given by
			\[
			f_{m-k}(B(w_{m-1})) = c(H^d_{w_{m-1}},k-1)\,.
			\]
			Meanwhile, the number of file placements of $m-k-1$ rooks on $B(w_{m})$ without a rook in the $m$-th column is equal to
			\[
			f_{m-k-1}(B(w_{m-1})) = c(H^d_{w_{m-1}},k)\,.
			\]
			An additional rook in the $m$-th column can then be added in $b_m$ ways to form a file placement of $m-k$ rooks on $B(w_m)$.
		\end{proof}
	\end{corollary}
	
	Our interpretation for the  $q$-normal order coefficients $\begin{bmatrix} w;q \\ k\end{bmatrix}$ involves the definitions that follow.
	
	\begin{definition} Let $G$ be a labeled digraph and let $\sigma$ be a cyclic partition of $G$.
		\begin{enumerate}
			\item The \emph{outdegree} of a vertex $i$, denoted by $\mathrm{deg}^{+}(i)$, is the number of vertices $j$ such that there is a directed edge from $i$ to $j$.
			\item The \emph{weight} of $i$ is defined as
			\[
			\mathrm{wt}^G_\sigma(i) =
			\begin{cases}
				i-1-\mathrm{deg}^{+}(i), & \mbox{if $i$ is a minimal element;}\\
				i-1-\mathrm{deg}^{+}(i) - \sigma^{(i)}(i),  & \mbox{otherwise}\,.
			\end{cases}
			\]
			\item The \emph{weight} of $\sigma$ is given by
			\[
			\mathrm{wt}^G(\sigma) = \sum_{i\in V(G)} \mathrm{wt}^G_\sigma(i)\,.
			\]
			\item The $q$-analogue of $c(G,k)$ is defined by
			\[
			c(G,k;q) = \sum_{\sigma} q^{\mathrm{wt}^G(\sigma)}\,.
			\]
			where the sum runs over all $k$-cyclic partitions $\sigma$ of $G$.
		\end{enumerate}
	\end{definition}

	\begin{theorem} \label{thm:firstq} Let $w$ be a Dyck word in $\{V,U\}$ subject to $UV=qVU+V$. Then $\begin{bmatrix} w;q \\ k \end{bmatrix}=c(H^d_w,k;q) = c(G^d_w,k;q)$.
		
		\begin{proof} Let $\sigma$ be a $k$-cyclic partition of $H^d_w$. We first show that for every $i$, $\mathrm{wt}^{H^d_w}_{\sigma}(i)$ is equal to the number of cells not lying above a rook in column $i$ of $\Phi(\sigma)$. This establishes that $\mathrm{wt}^{H_w}(\sigma) = \mathrm{wt}\left(\Phi(\sigma)\right)$ and together with Equation \eqref{eqn:celq}, proves that $\begin{bmatrix} w;q \\ k \end{bmatrix}=c(H^d_w,k;q)$.
			
			If $i$ is a minimal element of $\sigma$, then column $i$ of $\Phi(\sigma)$ does not contain a rook. The length of this row is equal to the number of vertices $j<i$ not adjacent to $i$. If $\mathrm{deg}^{+}(i)=d$, then the $d$ vertices adjacent to $i$ are $i-1,i-2,\ldots,i-d$. Equivalently, the vertices not adjacent to $i$ are $1,2,\ldots,i-d-1$, of which there are $i-1-d=\mathrm{wt}^{H^d_w}_{\sigma}(i)$.
			
			Next, if $i$ is not a minimal element of $\sigma$, then there is a rook in cell $(i,\sigma^{(i)}(i))$ of $\Phi(\sigma)$. As before, let $\mathrm{deg}^{+}(i)=d$. Since column $i$ has length $i-1-d$ by the argument in the previous paragraph, the number of cells below the rook is $i-1-d-\sigma^{(i)}(i)=\mathrm{wt}^{H^d_w}_{\sigma}(i)$. Thus, $\mathrm{wt}^{H^d_w}(\sigma)=\mathrm{wt}(\Phi(\sigma))$, as desired.
			
			Now, given a cyclic partition $\sigma$ of $H^d_w$ and its associated cyclic partition $\sigma'$ of $G^d_w$ as in the proof of Theorem~\ref{thm:firstmain}, the rook placements $\Phi(\sigma)$ and $\Phi(\sigma')$ have the same configuration up to a rearrangement of the rows, and thus, have the same weights, leading to $\mathrm{wt}^{H^d_w}{(\sigma)}=\mathrm{wt}^{G^d_w}{(\sigma')}$. This proves that $c(H^d_w,k;q) = c(G^d_w,k;q)$.
		\end{proof}	
	\end{theorem}	
	
	The following corollary gives a recurrence relation for $c(H^d_{w},k;q)$. Here, the $q$-analogue of a natural number $n$ is given by
	\[
	[n]_q =  1 + q + q^2 + \cdots + q^{n-1}\,.
	\]
	
	\begin{corollary} Let $w$ be a Dyck word in $\{V,U\}$ such that $|w|/2=n$, subject to $UV=qVU+V$. With $b_i$ and $w_i$ defined as in Corollary~\ref{cor:rec1}, for $k\leq m \leq n$,
		\[
		c(H^d_{w_m},k;q) = c(H^d_{w_{m-1}},k-1;q) + [b_{m}]_q \,c(H^d_{w_{m-1}},k;q)\,.
		\]
		\begin{proof}
			The proof is similar to that of Corollary~\ref{cor:rec1}. We simply note that the sum of the weights of all file placements of a rook in a column of length $b_m$ is equal to $[b_m]_q$.
		\end{proof}
	\end{corollary}
	
	\section{Linear partitions} \label{sec:line}
	
	The \emph{Lah numbers} $L(n,k)$ count the number of partitions of the set $\{1,2,\ldots,n\}$ into $k$ non-empty linearly-ordered subsets or lists. For example, the partitions of the set $\{1,2,3\}$ into two lists are $[1][2,3]$, $[1][3,2]$, $[2][1,3]$, $[2][3,1]$, $[3][1,2]$ and $[3][2,1]$ and thus, $L(3,2)=6$. 
	
	The Lah numbers satisfy the recurrence relation \[L(m,k)=L(m-1,k-1)+(m+k-1)L(m-1,k)\,,\] with initial conditions $L(m,0)=L(0,m)=\delta_{n,0}$, where $\delta_{\cdot,\cdot}$ denotes the Kronecker delta. These numbers also satisfy a simple explicit formula given by \[L(m,k)=\frac{m!}{k!}\binom{m-1}{k-1}\,.\]
	Goldman et al.~\cite[page 15]{Butetal} showed that $L(n,k)$ also counts the number of ways to place $n-k$ non-attacking rooks in a rectangular board with $n$ columns of height $n-1$. Other properties, combinatorial interpretations and $q$-analogues are discussed in \cite{Gon,Linetal}.
	
	In analogy with the interpretation of the Lah numbers, we define the following partition of a graph.
	
	\begin{definition}
		Let $G$ be a simple graph. A \emph{$k$-linear partition of $G$} is a partition of $V(G)$ into $k$ non-empty lists such that two elements $i$ and $j$ occur as consecutive elements in the same list only if vertices $i$ and $j$ are not adjacent in $G$. The {Lah number of a graph $G$}, denoted by $L(G,k)$, is the number of $k$-linear partitions of $G$.	
	\end{definition}
	
	For example, if $w=xxxDxDDxDDxxDD$, then a $3$-linear partition of the graph $H_w$ in Figure~\ref{fig:hw} is $[315][247][6]$ and a $3$-linear partition of the graph $G_w$ in Figure~\ref{fig:gw} is $[31][254][67]$. In addition, if $G$ is the empty graph with $n$ vertices, that is, $G=G_w=H_w$ with $w=(xD)^n$, then $L(n,k)=L(G,k)$. We can therefore think of $L(G,k)$ as a kind of generalization of the Lah numbers. 
	
	Let $w$ be a Dyck word in $\{x,D\}$ subject to $Dx=xD+1$ and $\mathcal X(w)$ be the word obtained by replacing every $x$ in $w$ with $x^2$. If $|w|/2=n$, we write the normally ordered form of $\mathcal X(w)$ by
	\[
	\mathcal X(w) = x^n \sum_{k=0}^n \left\lfloor
	\begin{matrix} 
		w \\ k
	\end{matrix}
	\right\rfloor	
	x^kD^k\,.
	\]
	
	A $q$-analog of the relation $Dx=xD+1$ is given by $Dx=qxD+1$. Two operators satisfying this relation are $x(f(x))=xf(x)$ and $D(x^m)=[m]_qx^{m-1}$. Given a Dyck word $w$ in $\{x,D\}$ with $|w|/2=n$, we write the normally ordered form of the word $\mathcal X(w)$ by
	\[
	\mathcal X(w) = x^n \sum_{k=0}^n \left\lfloor
	\begin{matrix} 
		w;q \\ k
	\end{matrix}
	\right\rfloor	
	x^kD^k\,.
	\]	
	
	Now, let $\psi$ be a placement of non-attacking rooks on a board $B$. The \emph{weight} of $\psi$, denoted by $\mathrm{wt}^*(\psi)$, is the number of cells not lying above or to the left of a rook. The $k$-th $q$-\emph{rook number} of $B$ is given by
	\[
	r_{k;q}(B) = \sum_{\psi} q^{\mathrm{wt}^*(\psi)}\,,
	\]
	where the sum runs over all placements $\psi$ of $k$ non-attacking rooks on $B$.	
	
	Using \cite[Identity 3]{Celetal}, we have
	\begin{align}
		\left\lfloor
		\begin{matrix} 
			w\\ k
		\end{matrix}
		\right\rfloor &= r_{n-k}(B(\mathcal X(w)))\,. \label{eq:lnkcela}\\
		\left\lfloor
		\begin{matrix} 
			w;q \\ k
		\end{matrix}
		\right\rfloor &= r_{n-k;q}(B(\mathcal X(w)))\,. \label{eq:lnkq}
	\end{align}

	Our main result for this section relates $\left\lfloor
	\begin{matrix} 
		w \\ k
	\end{matrix}
	\right\rfloor$ to the $k$-linear partitions of $H_w$ and $G_w$. First, we define the following map. Let $G$ be either $H_w$ or $G_w$ for some Dyck word $w$. Given a linear partition $\pi$ of $G$, define $\Psi(\pi)$ as the placement of non-attacking rooks on $B(\mathcal X(x))$ obtained as follows. For $i=1,2,\ldots,n$:
	\begin{enumerate}
		\item If $i$ is a minimal element, place no rook in column $i$.
		\item Let $\pi^{(i)}$ be the partial linear partition obtained by keeping only elements $j\leq i$. A position is \emph{allowable} if inserting $i$ into $\pi^{(i-1)}$ in that position does not result into $i$ occurring immediately before or immediately after a vertex adjacent to $i$. If $i$ is not a minimal element and $i$ is in the $t$-th position, that is, $\pi^{(i)}$ is obtained by inserting $i$ in the $t$-th allowable position of $\pi^{(i-1)}$ from the left, then place a rook in the $t$-th uncancelled cell from the bottom of column $i$.
	\end{enumerate} 
	In addition, for a simple graph $G$ with vertices labeled $1,2,\ldots,|V|$, let $\mathrm{deg}^<_G(i)$ be the number of vertices $j<i$ adjacent to vertex $i$.
	
	\begin{lemma} \label{lem:bijlnkx} Let $w$ be a Dyck word in $\{x,D\}$ with $|w|/2=n$. Then, the map $\Psi$ is a bijection between the $k$-linear partitions of $H_w$ (respectively, $G_w$) and the placements of $n-k$ non-attacking rooks on $B(\mathcal X(w))$. In other words,
		\begin{align*}
			L(H_w,k) &= r_{n-k}(B(\mathcal X(w)))\\
			L(G_w,k) &= r_{n-k}(B(\mathcal X(w)))\,.
		\end{align*}
		
		\begin{proof} We start by proving the first equation by showing that for a $k$-linear partition $\pi$ of $H_w$, the number of allowable positions to insert a non-minimal element $i$ in $\pi^{(i-1)}$ is equal to the number of uncancelled cells in column $i$ of $B(\mathcal X(w))$.
			
			Suppose that column $i$ of $B(w)$ has length $j$. Then, there is an edge between vertex $i$ and each of vertices $t=i-1,i-2,\ldots,j+1$. If $\mathrm{deg}^<_{H_w}(i)=d$, then there are $d=i-j-1$ such values of $t$ and $i$ cannot occur immediately before or immediately after each $t$. Suppose that $\pi^{(i-1)}$ contains $r$ minimal elements. Then, there are $(i-1+r)-2d=r-i+2j+1$ allowable positions to insert $i$ into $\pi^{(i-1)}$.
			
			Meanwhile, since column	$i$ of $B(w)$ has length $j$, column $i$ of $B(\mathcal X(w))$ has length $2j$. If $\pi^{(i-1)}$ has $r$ minimal elements, the first $i-1$ columns of $B(\mathcal X(w))$ contain $i-1-r$ non-attacking rooks. Therefore, column $i$ of $B(\mathcal X(w))$ has $2j-(i-1-r)=r-i+2j+1$ uncancelled cells, as desired.
			
			Lastly, if $\pi$ is a $k$-linear partition, then $\pi$ has $k$ minimal elements, so that $\Psi(\pi)$ contains $n-k$ rooks. This show that $\Psi$ is indeed a bijection, thereby establishing the first equation.
			
			Next, we prove the second equation. Observe by the previous argument that the number of allowable positions to insert an element $i$ in $\pi^{(i-1)}$ depends completely on $\mathrm{deg}^<_{H_w}(i)=d$ and  the number of minimal elements of $\pi^{(i-1)}$. This also holds for $G_w$. We now claim that for every $i$, $\mathrm{deg}^<_{H_w}(i)=\mathrm{deg}^<_{G_w}(i)$. It follows  that given a partial linear partition $\pi_1^{(i-1)}$ of $H_w$ and a partial linear partition $\pi_2^{(i-1)}$ of $G_w$ having the same minimal elements, the number of allowable positions of inserting $i$ into $\pi_1^{(i-1)}$ equals the number of allowable positions of inserting $i$ into $\pi_2^{(i-1)}$. Thus, if $\pi_1$ and $\pi_2$ have the same minimal elements and $i$ is in the same position in both $\pi_1^{(i)}$ and $\pi_2^{(i)}$ for every $i$, then $\Psi(\pi_1)=\Psi(\pi_2)$. For example, if $w=xxxDxDDxDDxxDD$, the linear partitions $\pi_1=[315][247][6]$ and $\pi_2=[31][254][67]$ correspond to the same rook placement shown in Figure~\ref{fig:hwgw}. The second equation now follows from the first.
			
			Finally, the claim is proved as follows. Observe that $\mathrm{deg}^<_{H_w}(i)$ is equal to the number of cells of column $i$ of $\mathcal B(H_w)$ and $\mathrm{deg}^<_{G_w}(i)$ is equal to the number of cells of column $i$ of $\mathcal B(G_w)$. By the first statement of Theorem~\ref{thm:alt}, $\mathrm{deg}^<_{H_w}(i)=\mathrm{deg}^<_{G_w}(i)$.
		\end{proof}
	\end{lemma}
	
	\begin{figure}[ht!] 
		\centering
		\begin{tikzpicture}[scale=0.6]
			\filldraw[color=gray!50, fill=gray!50] (0,6) rectangle (1,14);
			\filldraw[color=gray!50, fill=gray!50] (1,8) rectangle (3,14);
			\filldraw[color=gray!50, fill=gray!50] (3,10) rectangle (5,14);
			\filldraw[color=gray!50, fill=gray!50] (5,14) rectangle (7,14);
			
			\foreach \i in {1,...,7}{
				
				\node at (7.5-\i,14.5) {\scriptsize{$\i$}};
			}
			
			\foreach \i in {1,...,14}{
				\node at (-0.5,14.5-\i) {\scriptsize{$\i$}};
			}
			
			\foreach \i in {0,...,6}{
				\draw (\i,0) -- (\i,14);
				\draw (0,\i) -- (7,\i);
				\draw (0,7+\i) -- (7,7+\i);
			}

			\draw (7,0) -- (7,14);
			\draw (0,14) -- (7,14);
			
			\node at (4.5,10.5) {\Large{$\times$}};
			\node at (3.5,13.5) {\Large{$\times$}};
			\node at (2.5,11.5) {\Large{$\times$}};
			\node at (0.5,12.5) {\Large{$\times$}};
		\end{tikzpicture}
		
		\caption{The rook placement associated with both $\pi_1=[315][247][6]$ and $\pi_2=[31][254][67]$ as $3$-linear partitions of $H_w$ and $G_w$, respectively, with $w=xxxDxDDxDDxxDD$.} \label{fig:hwgw}
		
	\end{figure}
	
	The next theorem states the main result for this section, which is a combinatorial interpretation for the normal order coefficients $\left\lfloor
	\begin{matrix} 
		w \\ k
	\end{matrix}
	\right\rfloor$ of $\mathcal X(w)$ in terms of linear partitions of $H_w$ and $G_w$.
	
	\begin{theorem} \label{thm:lnkmain}Let $w$ be a Dyck word in $\{x,D\}$ subject to $Dx=xD+1$. Then, $\left\lfloor
		\begin{matrix} 
			w \\ k
		\end{matrix}
		\right\rfloor = L(H_w,k) = L(G_w,k)$
		\begin{proof}
			The result follows from Lemma~\ref{lem:bijlnkx} and Equation~\eqref{eq:lnkcela}.
		\end{proof}
	\end{theorem}
	
	The next corollary gives a recurrence relation for $L(H_{w_n},k)$. As with the previous recurrence relations, the proof is similar to that of Corollary~\ref{cor:rec1}. Observe that if $w=(xD)^n$, then $b_m=m-1$ and we recover the recurrence relation for the classical Lah numbers.
	
	\begin{corollary} \label{cor:qlnka} Let $w$ be a Dyck word in $\{x,D\}$ such that $|w|/2=n$, subject to $Dx=xD+1$. Suppose that the column lengths of $B(w)$ are $b_1,b_2,\ldots,b_n$. Let $w_i$ be a subword of $w$ such that the board $B(w_i)$ consists of columns $1,2,\ldots,i$ of $B(w)$. Then for $k\leq m\leq n$, 
		\[
		L(H^d_{w_m},k) = L(H^d_{w_{m-1}},k-1) + \left[2b_{m}-(m-k-1)\right] L(H^d_{w_{m-1}},k)\,.
		\]	
	\end{corollary}
	
	Before stating our combinatorial interpretation for $\left\lfloor
	\begin{matrix} 
		w;q \\ k
	\end{matrix}
	\right\rfloor$, we will need the following definitions.
	
	\begin{definition} Let $\pi$ be a linear partition of a simple graph $G$ with vertices labeled $1,2,\ldots,|G|$.
		\begin{enumerate}	
			\item The \emph{weight} of a vertex $i$ is defined by
			\[
			\mathrm{wt}_{\pi}^{*G}(i)  =
			\begin{cases}
				2\bigl(i-1-\mathrm{deg}^{<}_{G}(i)\bigr) - |\{j \mbox{ non-minimal},j<i\}|, & \mbox{$i$ is a minimal element}\,;\\
				t-1, $\mbox{ $i$ is in position $t$ in $\pi^{(i)}$}$, &\mbox{otherwise}\,.			
			\end{cases}
			\]	
			\item The \emph{weight} of $\pi$ is defined by
			\[
			\mathrm{wt}^{*G}(\pi)=\sum_{i\in V(G)} \mathrm{wt}_{\pi}^{*G}(i)\,.
			\]
			\item The $q$-analogue of $L(G,k)$ is given by
			\[
			L(G,k;q) = \sum_{\pi} q^{\mathrm{wt}^{*G}(\pi)}\,,
			\]
			where the sum runs over all $k$-linear $\pi$ partitions of $G$.
		\end{enumerate}
	\end{definition}

	\begin{theorem} Let $w$ be a Dyck word in $\{x,D\}$ subject to $Dx=qxD+1$. Then $\left\lfloor
		\begin{matrix} 
			w;q \\ k
		\end{matrix}
		\right\rfloor = L(H_w,k;q) = L(G_w,k;q)$.
		
		\begin{proof} To prove $\left\lfloor
			\begin{matrix} 
				w;q \\ k
			\end{matrix}
			\right\rfloor = r_{n-k;q}(B(\mathcal X(w)))$, we show that given a linear partition $\pi$ of $H_w$, the number of uncancelled cells of column $i$ of $\Psi(\pi)$ is equal to $\mathrm{wt}^{*H_w}_\pi(i)$. If $i$ is a minimal element, then column $i$ of $B(\mathcal X(w))$ has $2\bigl(i-1-\mathrm{deg}^{<}_{G}(i)\bigr)$ cells. The number of rooks to the right of this column is equal to the number of non-minimal elements $j<i$. Meanwhile, if $i$ is not a minimal element that is in position $t$ in $\pi^{(i)}$, then there is a rook in the $t$-th uncancelled cell from the bottom of column $i$. Below such rook are $t-1$ uncanceled cells. In combination with Equation~\ref{eq:lnkq}, this proves $\left\lfloor
			\begin{matrix} 
				w;q \\ k
			\end{matrix}
			\right\rfloor = L(H_w,k;q)$.
			
			Lemma~\ref{lem:bijlnkx} shows that given a $k$-linear partition $\pi_1$ of $H_w$, there exists a unique $k$-linear partition $\pi_2$ of $G_w$ such that $\Psi(\pi_1)=\Psi(\pi_2)$. This, together with $\left\lfloor
			\begin{matrix} 
				w;q \\ k
			\end{matrix}
			\right\rfloor = L(H_w,k;q)$, establishes $L(H_w,k;q) = L(G_w,k;q)$
		\end{proof}
	\end{theorem}
	
	The corollary that follows gives a recurrence relation for $L(H^d_{w_m},k)$. 
	
	\begin{corollary} Let $w$ be a Dyck word in $\{x,D\}$ such that $|w|/2=n$, subject to $Dx=qxD+1$. With $b_i$ and $w_i$ defined as in Corollary~\ref{cor:qlnka}, for $k\leq m \leq n$, 
		\[
		L(H_{w_m},k) = L(H_{w_{m-1}},k-1) + \left[2b_{m}-(m-k-1)\right]_q L(H_{w_{m-1}},k)\,.
		\]	
	\end{corollary}
	
	We note that setting $w=(xD)^n$ yields a $q$-analogue of the Lah numbers with recurrence relation
	\[
	L_q[n,k] = L_q[n-1,k-1] + ([n+k-1]_q) L_q[n-1,k]\,.
	\]
	This $q$-analogue is slightly different from the one studied by Lindsay et al.~\cite{Linetal} and Gonzales~\cite{Gon}, where the $q$-Lah numbers satisfy
	\[
	L^*_q[n,k] = L^*_q[n-1,k-1] + ([n-1]_q+[k_q]) L^*_q[n-1,k]\,.
	\]
	
	\section{Some Remarks} \label{sec:sumrem}
	
	\begin{enumerate}
		\item The results in this paper assume that $w$ is a Dyck word. Using the construction used by Engbers et al.~\cite[Definition 2.4]{Engetal}, we can generalize our results to arbitrary words. This is done by rewriting $w$ into $w'=x^awD^b$ (or $V^awU^b$), where $a$ and $b$ are the smallest integers such that $w'$ is a Dyck word. We then have
		\begin{align*}
			\begin{bmatrix} 
				w\\ k
			\end{bmatrix} &= c(H^d_{w'},k+b) \\	
			\left\lfloor
			\begin{matrix} 
				w\\ k
			\end{matrix}
			\right\rfloor &= L(H_{w'},k+b)\,.	
		\end{align*}
		The same modification applies to the $q$-analogues.
		\item A different graphical Stirling number of the first kind $c^*(G,k)$ was studied by Barghi~\cite{Bar}, where $c^*(G,k)$ counts the number of permutations $\sigma$ of $V(G)$ into $k$ cycles such that $v\sigma(v)$ is an edge in the graph obtained by adding a loop at each vertex of $G$. This number seems to be less restrictive than $c(G,k)$, and it would be interesting to find an interpretation for $c^*(G,k)$ as normal order coefficients. 
		\item It follows from Theorem~\ref{thm:alt} that other graphs $G$ for which $\displaystyle {w\brace k}=S(G,k)$ can be constructed by rearranging the rows of $B(w)$ as long as the cells remain above the northeast diagonal of the lattice. In addition, the rook equivalence of boards having different column and row lengths~\cite[Section 1.6]{Butetal} leads to still other graphs. It would be interesting to determine the rules for constructing such graphs from a Dyck word $w$ and derive necessary and sufficient conditions to determine graphs $G_1$ and $G_2$ such that $S(G_1,k)=S(G_2,k)$. How many isomorphism classes satisfy $S(G_1,k)=S(G_2,k)$ are there given $w$?
	\end{enumerate}

\end{document}